\documentclass[reqno]{amsart}
\usepackage[utf8]{inputenc}

\usepackage{caption}
\usepackage{amssymb}
\usepackage{graphicx}
\usepackage{latexsym}
\usepackage{stmaryrd}
\usepackage{enumitem}
\usepackage[bookmarks]{hyperref}
\usepackage{multirow}
\usepackage{xcolor}
\usepackage{relsize}
\usepackage{comment}
\usepackage{xifthen}
\usepackage{mathrsfs}
\usepackage{svg}
\usepackage{mathtools}
\usepackage{tikz}
\usepackage{wrapfig}
\usetikzlibrary{arrows}
\usetikzlibrary{positioning}
\colorlet{darkishRed}{red!80!black}
\colorlet{darkishBlue}{blue!60!black}
\colorlet{darkishGreen}{green!60!black}
\hypersetup{
    draft = false,
    bookmarksopen=true,
    colorlinks,
    linkcolor={red!60!black},
    citecolor={green!60!black},
    urlcolor={blue!60!black}
}
\usepackage{geometry}
\usetikzlibrary{arrows}

\newcommand{\nao}[1][]{%
\ifthenelse{\equal{#1}{}}{\trianglelefteq_T}{\trianglelefteq_{T_{#1}}\!}%
}

\def\lowfwd #1#2#3{{\mathop{\kern0pt #1}\limits^{\kern#2pt\raise.#3ex
\vbox to 0pt{\hbox{$\scriptscriptstyle\rightarrow$}\vss}}}}
\def\lowbkwd #1#2#3{{\mathop{\kern0pt #1}\limits^{\kern#2pt\raise.#3ex
\vbox to 0pt{\hbox{$\scriptscriptstyle\leftarrow$}\vss}}}}

\def\ve{\kern-1.5pt\lowfwd e{1.5}2\kern-1pt}
\def\ev{\kern-1pt\lowbkwd e{0.5}2\kern-1pt}
\def\vf{\kern-2pt\lowfwd f{2.5}2\kern-1pt}

\renewcommand{\subset}{\subseteq}

\newcommand{ \N } { \mathbb{N} }

\makeatletter

\def\calCommandfactory#1{%
   \expandafter\def\csname c#1\endcsname{\mathcal{#1}}}
\def\frakCommandfactory#1{%
   \expandafter\def\csname frak#1\endcsname{\mathfrak{#1}}}

\newcounter{ctr}
\loop
  \stepcounter{ctr}
  \edef\X{\@Alph\c@ctr}
  \expandafter\calCommandfactory\X
  \expandafter\frakCommandfactory\X
  \edef\Y{\@alph\c@ctr}
  \expandafter\frakCommandfactory\Y
\ifnum\thectr<26
\repeat

\newtheorem*{HalinThm}{Halin's grid theorem}
\newtheorem{mainresult}{Theorem}

\theoremstyle{definition}

\theoremstyle{remark}

\title[A strengthening of Halin's grid theorem]{A strengthening of Halin's grid theorem\\}

\author{Jan Kurkofka, Ruben Melcher, and Max Pitz}
\address{Universit\"at Hamburg, Department of Mathematics, Bundesstrasse 55 (Geomatikum), 20146 Hamburg, Germany}
\email{\{jan.kurkofka, ruben.melcher, max.pitz\}@uni-hamburg.de}

\@namedef{subjclassname@2020}{\textup{2020} Mathematics Subject Classification}

\begin{document}

\begin{abstract}
We show that for every infinite collection $\cR$ of disjoint equivalent rays in a graph $G$ there is a subdivision of the hexagonal half-grid in $G$ such that all its vertical rays belong to $\cR$. This result strengthens Halin's grid theorem by giving control over which specific set of rays is used, while its proof is significantly shorter.
\end{abstract}
\vspace*{-2cm}
\maketitle

\vspace*{-.7cm}

\section{Introduction}
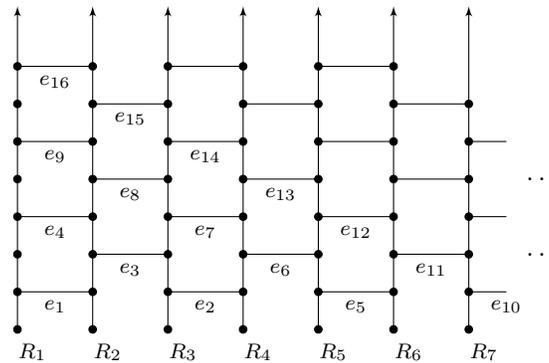
\begin{figure}[ht]
    \begin{tikzpicture}
    \tikzset{edge/.style = {->,> = latex'}}
    
    \foreach \x in {-4,-2,0} {
    \draw[edge] (\x,0) to (\x,4.3);
    
     \draw[edge] (\x +1,0) to (\x +1 ,4.3);
     
     \foreach \y in {0.5,1.5,2.5,3.5} {
      \draw[fill,black] (\x,\y) circle (.05);
       \draw[fill,black] (\x+1,\y) circle (.05);
         \draw (\x,\y) to (\x+1,\y);
     }
    
        }

     \foreach \x in {-3,  -1 } {
         \foreach \y in {1,2,3}{

         \draw[fill,black] (\x,\y) circle (.05);
         
         \draw[fill,black] (\x +1 ,\y) circle (.05);
         
         \draw (\x,\y) to (\x+1,\y);
        }
         }
         

\draw[edge] (2,0) to (2,4.3);
   \foreach \y in {1,2,3}{

         \draw[fill,black] (1,\y) circle (.05);
          \draw[fill,black] (2,\y) circle (.05);
           \draw (1,\y) to (2,\y);
            \draw[fill,black] (2,\y-0.5) circle (.05);
             \draw (2,\y-0.5) to (2.5,\y-0.5);
         }

    \node at (3,2) {$\dots$};
     \node at (3,1) {$\dots$};

    \foreach \x in {-4,...,2}{
     \draw[fill,black] (\x,0) circle (.05);
     }
     
\foreach \x in {1,...,7}{
    \node at (\x-4.8,-0.3) {\footnotesize{$R_{\x}$}};
     }
    
    \node at (-3.5, 0.3 ) {\footnotesize{$e_1$}};
     \node at (-3.5, 1.3 ) {\footnotesize{$e_4$}};
      \node at (-3.5, 2.3 ) {\footnotesize{$e_9$}};
       \node at (-3.5, 3.3 ) {\footnotesize{$e_{16}$}};
       
    
    \node at (-2.5, 0.8 ) {\footnotesize{$e_3$}};
     \node at (-2.5, 1.8 ) {\footnotesize{$e_8$}};
      \node at (-2.5, 2.8 ) {\footnotesize{$e_{15}$}};
      
    
    \node at (-1.5, 0.3 ) {\footnotesize{$e_2$}};
     \node at (-1.5, 1.3 ) {\footnotesize{$e_7$}};
      \node at (-1.5, 2.3 ) {\footnotesize{$e_{14}$}};
    
    \node at (-.5, 0.8 ) {\footnotesize{$e_6$}};
     \node at (-.5, 1.8 ) {\footnotesize{$e_{13}$}};

    
    \node at (.5, 0.3 ) {\footnotesize{$e_5$}};
     \node at (.5, 1.3 ) {\footnotesize{$e_{12}$}};
     
\node at (1.5, 0.8 ) {\footnotesize{$e_{11}$}};
     \node at (2.5, .3 ) {\footnotesize{$e_{10}$}};


\foreach \y in {1,2,3}{

         \draw[fill,black] (-4,\y) circle (.05);
         }

\end{tikzpicture}
\caption{The hexagonal half-grid with vertical rays $R_i$.}
\label{fig:my_label}
\end{figure}

\noindent  An \emph{end} of a graph $G$ is an equivalence class of rays, where two rays of $G$ are \emph{equivalent} if there are infinitely many vertex-disjoint paths between them in~$G$. The \emph{degree} $\deg(\omega) \in \N \cup \{\infty\}$
of an end $\omega$ of $G$ is the maximum size of a collection of pairwise disjoint rays in~$\omega$, see~Halin~\cite{H65}. Ends of infinite degree are also called \emph{thick}. 
The \emph{half-grid}, the graph on $\N^2$ in which two vertices $(n,m)$ and $(n',m')$ are adjacent if and only if $|n-n'|+|m-m'| = 1$, and its sibling the \emph{hexagonal half-grid}, where one deletes every other rung from the half-grid as shown in Figure~\ref{fig:my_label}, are examples of graphs which have only one end, which is thick.

\noindent One of the cornerstones of infinite graph theory, \emph{Halin's grid theorem}~\cite{H65},  says that grid-like graphs are the prototypes for ends of infinite degree: 

\begin{HalinThm}
 Every graph with an end of infinite degree contains a subdivision of the hexagonal half-grid whose rays belong to that end. 
\end{HalinThm}

Halin's theorem is a precursor of the work by Robertson, Seymour and Thomas on excluding infinite grid or clique minors~\cite{robertson1995excluding} and has further influenced research in \cite{bowler2018ubiquity,gkmp,halin1975problem,heuer2017excluding}. It is curious, however, that Halin's theorem does not mention any specific ray families, that is if one chooses a specific infinite collection $\cR$ of  disjoint rays witnessing that the end is thick, then neither the assertion of Halin's theorem nor its available proofs by Halin~\cite[Satz~4]{H65} and by Diestel~\cite{dlestel2004short,Bible} make any assertion on how the resulting subdivided hexagonal half-grid relates to the collection of rays $\cR$ one started with. Furthermore, in recent work of ours on an extension of Halin's grid theorem to higher cardinals \cite{gkmp}, it became quite important to achieve more control of specific uncountable ray families, and the question arose whether this can be done also in the countable case. And indeed, the main result of this note is that this is in fact possible:

\begin{mainresult}\label{mainresult}
For every infinite collection $\cR$ of disjoint equivalent rays in a graph $G$ there is a subdivision of the hexagonal half-grid in $G$ such that all its vertical rays belong to $\cR$.
\end{mainresult}

The  known proofs of Halin's  grid theorem are rather involved and require an elaborate recursive construction that runs close to five pages in Diestel's textbook. Our stronger result in Theorem~\ref{mainresult} requires a different  approach -- which coincidentally provides a much shorter proof of Halin's original grid theorem.

\section{The proof}

\noindent Suppose we are handed a countably infinite collection $\cR$ of disjoint equivalent rays in a graph $G$. By routine arguments, we fix, for the remainder of this paper, a ray $S$ in $G$ that meets each ray in $\cR$ infinitely often. 
An \emph{$\cR$-segment (of $S$)} is any maximal subpath of $S$ which is internally disjoint from $\bigcup\cR$. We say that an $\cR$-segment is  \emph{between} two rays $R_1,R_2 \in \cR$ if it has its endpoints on $R_1$ and $R_2$ respectively. Let $M(\cR)$ denote the auxiliary multigraph with vertex  set $\cR$ where the multiplicity of an edge $R_1 R_2$ is equal to the number of $\cR$-segments between $R_1$~and~$R_2$. Finally, let $M_\infty(\cR)$ denote the spanning subgraph of $M(\cR)$ obtained by removing all edges of finite multiplicity.

\begin{center}
    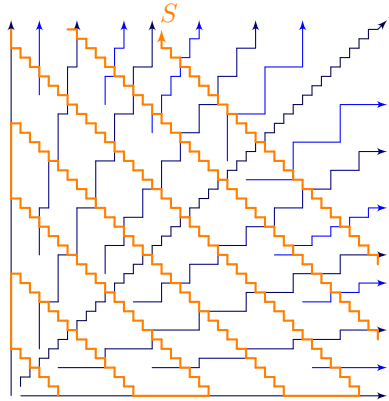
\begin{figure}[h]
        \centering
\begin{tikzpicture}
\tikzset{edge/.style = {->,> = latex'}}
\tikzset{line/.style = {-,line cap=round}}
    
    \def\drawLine#1#2#3#4#5#6{
    \edef\x0{#1}
    \edef\y0{#2}
    \edef\dx{#3}
    \edef\dy{#4}
    \edef\n{#5}
    \edef\col{#6}
    \pgfmathparse{int(\n-1)}
    \edef\m{\pgfmathresult}
    \foreach \t in {0,...,\m}{
        \draw[line,color=black!\col!blue] ({\x0+\t*\dx},{\y0+\t*\dy}) to ({\x0+(\t+1)*\dx},{\y0+\t*\dy});
        \draw[line,color=black!\col!blue] ({\x0+(\t+1)*\dx},{\y0+\t*\dy}) to ({\x0+(\t+1)*\dx},{\y0+(\t+1)*\dy});
    }
    \draw[edge,color=black!\col!blue] ({\x0+\n*\dx},{\y0+\n*\dy}) to ({\x0+(\n+1)*\dx},{\y0+\n*\dy});
    }
    
    \def\drawLineI#1#2#3#4#5#6{
    \edef\x0{#2}
    \edef\y0{#1}
    \edef\dx{#4}
    \edef\dy{#3}
    \edef\n{#5}
    \edef\col{#6}
    \pgfmathparse{int(\n-1)}
    \edef\m{\pgfmathresult}
    \foreach \t in {0,...,\m}{
        \draw[line,color=black!\col!blue] ({\x0+\t*\dx},{\y0+\t*\dy}) to ({\x0+\t*\dx},{\y0+(\t+1)*\dy});
        \draw[line,color=black!\col!blue] ({\x0+\t*\dx},{\y0+(\t+1)*\dy}) to ({\x0+(\t+1)*\dx},{\y0+(\t+1)*\dy});
    }
    \draw[edge,color=black!\col!blue] ({\x0+\n*\dx},{\y0+\n*\dy}) to ({\x0+\n*\dx},{\y0+(\n+1)*\dy});
    }
    
    \def\drawLineM#1#2#3#4#5#6{
    \edef\x0{#2}
    \edef\y0{#1}
    \edef\dx{#4}
    \edef\dy{#3}
    \edef\n{#5}
    \edef\col{#6}
    \pgfmathparse{int(\n-1)}
    \edef\m{\pgfmathresult}
    \foreach \t in {0,...,\m}{
        \draw[line,color=black!\col!blue] ({\x0+\t*\dx},{\y0+\t*\dy}) to ({\x0+\t*\dx},{\y0+(\t+1)*\dy});
        \draw[line,color=black!\col!blue] ({\x0+\t*\dx},{\y0+(\t+1)*\dy}) to ({\x0+(\t+1)*\dx},{\y0+(\t+1)*\dy});
    }
    }
    
    \def\drawLineS#1#2#3#4#5{
    \edef\x0{#1}
    \edef\y0{#2}
    \edef\dx{#3}
    \edef\dy{#4}
    \edef\n{#5}
    \foreach \t in {0,...,\n}{
        \draw[line,orange,thick] ({\x0+\t*\dx},{\y0+\t*\dy}) to ({\x0+(\t+1)*\dx},{\y0+\t*\dy});
        \draw[line,orange,thick] ({\x0+(\t+1)*\dx},{\y0+\t*\dy}) to ({\x0+(\t+1)*\dx},{\y0+(\t+1)*\dy});
    }
    }
    
    \draw[edge,color=black!70!blue] (0,0) to (0,5);
    \draw[edge,color=black!70!blue] (.125,0) to (5,0);
    \drawLineM{.125}{.125}{.125}{.125}{38}{70}
    \draw[edge,color=black!70!blue] (4.875,4.875) -- (5,5);
    
    \drawLine{.625}{.375}{.625}{.25}{6}{60}
    \drawLineI{.625}{.375}{.625}{.25}{6}{60}

    \drawLine{1.875}{.375}{.625}{.125}{4}{50}
    \drawLineI{1.875}{.375}{.625}{.125}{4}{50}
    
    \drawLine{1.625}{1.25}{.375}{.25}{8}{50}
    \drawLineI{1.625}{1.25}{.375}{.25}{8}{50}
    
    \drawLine{3.125}{2.875}{.625}{.5}{2}{0}
    \drawLineI{3.125}{2.875}{.625}{.5}{2}{0}
    
    \draw[edge,color=black!0!blue] (4,.375) -- (5,.375);
    \draw[edge,color=black!0!blue] (.375,4) -- (.375,5);
    
    \drawLine{3.875}{1.25}{.375}{.125}{2}{0}
    \drawLineI{3.875}{1.25}{.375}{.125}{2}{0}
    
    \drawLine{3.5}{1.875}{.25}{.125}{5}{0}
    \drawLineI{3.5}{1.875}{.25}{.125}{5}{0}
    
    \drawLineS{0}{.625}{.125}{-.125}{4}
    \draw[line,orange,thick] (0,.625)--(0,1.625);
    \drawLineS{0}{1.625}{.125}{-.125}{12}
    \draw[line,orange,thick] (1.625,0)--(2.625,0);
    \drawLineS{0}{2.625}{.125}{-.125}{20}
    \draw[line,orange,thick] (0,2.625)--(0,3.625);
    \drawLineS{0}{3.625}{.125}{-.125}{28}
    \draw[line,orange,thick] (3.625,0)--(4.625,0);
    \drawLineS{0}{4.625}{.125}{-.125}{36}
    \draw[line,orange,thick] (0,4.625)--(0,4.875);
    \drawLineS{.75}{4.875}{.125}{-.125}{32}
    \drawLineS{2}{4.625}{.125}{-.125}{22}
    \draw[edge,orange,thick] (2,4.625)--(2,4.875);
    \node[orange] at (2.1,5.1) {$S$};
    
\end{tikzpicture}
        \caption{A configuration of rays in the half-grid with edge-less $M_\infty(\cR)$.}
        \label{fig:my_label2}
    \end{figure}
\end{center}

Since $S$ meets every ray in $\cR$ infinitely often, it follows that $M(\cR)$ is infinitely edge-connected. 
Further, if $M_\infty(\cR)$ has a component with infinitely many vertices, then Theorem~\ref{mainresult} follows at once: 
In this case, $M_\infty(\cR)$ either contains a ray $R_1,R_2,\ldots$ or an infinite star with centre $R$ and leaves $R_1,R_2,\ldots$. In the ray case, one recursively selects sufficiently late $\cR$-segments to represent subdivided edges $e_1,e_2,e_3,\ldots$ of the hexagonal half-grid in the order indicated in Figure~\ref{fig:my_label}; in the star case, edges between $R_i$ and $R_j$ are represented by two sufficiently late $\cR$-segments between $R_i, R_j$ and  $R$ together with the subpath on $R$ connecting the endpoints of those segments.

However, $M_\infty(\cR)$ might have no edges at all: Consider for example a collection of radial rays in the  half-grid such that between any two rays there lies a third, see Figure~\ref{fig:my_label2}. Still, by moving to an infinite subcollection $\cR' \subseteq \cR$ and considering auxiliary multigraphs $M(\cR')$ and $M_\infty(\cR')$ instead, the connectivity properties of $M_\infty(\cR')$ might improve. Indeed, the auxiliary multigraphs for $\cR'$ remain well-defined as the same $S$ still meets every ray in $\cR'$ infinitely often.
Note, however, that $\cR$-segments of $S$ may now be properly contained in $\cR'$-segments of $S$. Our preceding discussion can be summarized as:
\begin{enumerate}
    \item \label{item_1} The auxiliary multigraph $M(\cR')$ is infinitely edge-connected for any infinite $\cR' \subset \cR$.
    \item \label{item_2} If for some $\cR' \subset \cR$ the auxiliary multigraph $M_\infty(\cR')$ has an infinite component, there is a subdivision of the hexagonal half-grid in $G$ such that all its vertical rays belong to~$\cR'$.
\end{enumerate}
Our next observation provides a  sufficient condition for $M_\infty(\cR')$ to have an infinite component. Recall that the degree of a vertex in a multigraph denotes the number of its neighbours.
\begin{enumerate}[resume]
    \item \label{item_3} If $\cR' \subseteq \cR$ is infinite such that $M(\cR')$ has only finitely many vertices of infinite degree, then $M_\infty(\cR')$ has an infinite component.
\end{enumerate}
Indeed, suppose for a contradiction that all components of $M_\infty(\cR')$ are finite. 
Then there is also a finite component $C$ of $M_\infty(\cR')$ that contains none of the finitely many vertices that have infinite degree in $M(\cR')$. 
Since $M(\cR')$ is infinitely edge-connected by (\ref{item_1}), there are infinitely many edges in $M(\cR')$ from $C$ to its complement. And since $C$ consists of vertices of finite degree only, the neighbourhood of $C$ in $M(\cR')$ is finite. Thus, there is a vertex in $C$ that sends infinitely many edges to some vertex outside of~$C$, contradicting the choice of $C$. This establishes (\ref{item_3}).

\medskip

The idea of the proof of Theorem~\ref{mainresult} is now as follows: If all vertices of $M(\cR)$ have finite degree, then $M_\infty(\cR)$ has an infinite component by (\ref{item_3}) and we are done by (\ref{item_2}). Otherwise, there is a ray $R_1 \in \cR$ that has infinitely many neighbours $N(R_1)$ in $M(\cR)$, and we may restrict our collection of rays to $\cR_1 := \{R_1\} \cup N(R_1)$. 
Next, if all but finitely many rays of $\cR_1$ have finite degree in $M(\cR_1)$, then $M_\infty(\cR_1)$ has an infinite component by (\ref{item_3}) and we are again done by (\ref{item_2}). Thus, we may pick a second ray $R_2$ in $\cR_1$ distinct from $R_1$ such that $N(R_2)$ is infinite in $M(\cR_1)$, and restrict our collection of rays to $\cR_2 := \{R_1,R_2\} \cup N(R_2)$. 
Repeating this step as often as possible gives rise to a sequence of rays $R_1,R_2,R_3,\ldots$. If this procedure ever stops because there are no more vertices of infinite degree to choose, then we are done by (\ref{item_3})  and (\ref{item_2}). Thus, the question becomes what to do when this procedure does not terminate. 

Informally, the solution is to modify our construction so that besides the first $n$ rays  $R_1,\ldots,R_n$ we will also have chosen suitable paths $P_1,\ldots,P_{n-2}$ between them representing the subdivided edges $e_1,\ldots,e_{n-2}$ in the copy of the hexagonal half-grid from Figure~\ref{fig:my_label}.\footnote{The index shift just has the purpose that when choosing a path for $e_2$ we have already selected $R_3$ and $R_4$.} Then, in the case where our procedure never stops, the chosen rays $R_1,R_2,\ldots$ become the vertical rays of a hexagonal half-grid where the subdivided paths corresponding to an edge $e_i$ are given by the path $P_i$.

Formally, suppose that at step $n$ we have chosen $n$ distinct rays $R_1,\ldots,R_n$ from $\cR$ and an infinite subcollection $\cR_n \subset \cR$ containing all chosen $R_i$ such that in $M(\cR_n)$ every $R_i$ is adjacent to all rays in $\cR_n \setminus \{R_1,\ldots,R_n\}$. Further, suppose that we have chosen $n-2$ disjoint paths $P_1,\ldots,P_{n-2}$ internally disjoint from $\bigcup \cR_n$, such that each $P_i$ connects the same two rays from $\{R_1,\ldots,R_n\}$ as $e_i$ in Figure~\ref{fig:my_label}, in a way such that whenever two paths $P_i,P_j$ with $i<j$ have endvertices on the same ray $R_k$, then the endvertex of $P_i$ comes before the endvertex of $P_j$ on $R_k$. 

Now if all but finitely many rays in $\cR_n$ have finite degree in $M(\cR_n)$, then we are done by (\ref{item_3})  and (\ref{item_2}). Hence, we may assume that there is a ray $R_{n+1}$ in $\cR_n \setminus \{R_1,\ldots,R_n\}$ that has infinitely many neighbours $N(R_{n+1})$ in $M(\cR_n)$. Let $\cR'_{n+1} :=   \{R_1, \ldots, R_{n+1}\} \cup N(R_{n+1})$, and note that in $M(\cR'_{n+1})$, every $R_i$ for $i=1,\ldots,n+1$ is adjacent to all other rays in $\cR'_{n+1} $. Now let $i$ and $j$ denote the indices of the rays in Figure~\ref{fig:my_label} containing the endvertices of the edge $e_{n-1}$.
Note that $i,j \leq n+1$. Since both $R_i$ and $R_j$ have infinitely many common neighbours in $M(\cR'_{n+1})$, we also find a common neighbour $Q_{n-1}$ in $\cR'_{n+1} \setminus \{R_1,\ldots,R_{n+1}\}$ such that the corresponding $\cR'_{n+1}$-segments of $S$ between $R_i,R_j$ and $Q_{n-1}$ are disjoint from all earlier paths $P_1,\ldots,P_{n-2}$ and also have their endvertices on $R_i,R_j$ later than the endvertices of any previous path $P_1,\ldots,P_{n-2}$. Then we may pick a new path $P_{n-1}$ consisting of both these $\cR'_{n+1}$-segments of $S$ between $R_i,R_j$ and $Q_{n-1}$ together with a suitable subpath of $Q_{n-1}$. Finally, set $\cR_{n+1}:= \cR'_{n+1} \setminus \{Q_{n-1}\}$. This completes the induction step, and the proof is complete. \hfill \qed

\, \newline
We remark that only in the case where our procedure stops and~(\ref{item_2}) yields a ray one can just build a grid, in all other cases one can build a clique of rays, i.e.\ one finds an infinite $ \cR' \subseteq \cR$ and a family of internally disjoint $\bigcup\cR'$-paths witnessing that any two rays in $\cR'$ are equivalent.

\bibliographystyle{plain}
\bibliography{reference}

\end{document}